\documentclass[11pt,epsfig,a4paper]{article}
\usepackage{amssymb}
\usepackage{amsmath}
\usepackage{amsthm}
\usepackage{epsfig}
\usepackage{longtable}
\newtheoremstyle{thmm}{1.5ex plus 1ex minus .2ex}{1.5ex plus 1ex minus
.2ex}{\rmfamily}{}{\bfseries}{}{1em}{}
\theoremstyle{thmm}
\newtheorem{theorem}{Theorem}[section]
\newtheorem{lemma}{Lemma}[section]

\renewenvironment{proof}[1][Proof]{\noindent\textit{#1. } }{\hfill$\square$}
\allowdisplaybreaks \textwidth  6in \textheight 9in \topmargin
-0.1in \oddsidemargin 0.2in \evensidemargin 0.0in

\newcommand{\nn}{\nonumber}
\def \endproof{\vrule height8pt width 5pt depth 0pt}

\def\d{\delta}

\def\R{\mathbb{R}}

\def\d{\,{\rm d}}

\def\u{{\bf u}}

\begin{document}
\title{\bf
Unconditional convergence and optimal error estimates of a
Galerkin-mixed FEM for incompressible miscible flow in porous media}
\author{
Buyang Li\footnote{Department of Mathematics,
Nanjing University, Nanjing, China.
 {\tt buyangli@nju.edu.cn}
}
~  and~ Weiwei Sun\footnote{Department of Mathematics,
City University of Hong Kong, Hong Kong.  {\tt maweiw@math.cityu.edu.hk}
\newline\indent~
The work of the
authors was supported in part by a grant from the Research Grants
Council of the Hong Kong Special Administrative Region, China
(Project No. CityU 102005)
\indent~ }
}
\date{}

\maketitle

\begin{abstract}
In this paper, we study the unconditional convergence and error
estimates of a Galerkin-mixed FEM with the linearized semi-implicit
Euler time-discrete scheme for the equations of incompressible
miscible flow in porous media. We prove that the optimal $L^2$ error
estimates hold without any time-step (convergence) conditions, while
all previous works require certain time-step restrictions.
Our theoretical results provide a new understanding on commonly-used
linearized schemes. The proof is based on a splitting of the error
into two parts: the error from the time discretization of
the PDEs and the error from the finite element discretization of
corresponding time-discrete PDEs. The approach used in this paper can be
applied to more general nonlinear parabolic systems and many
other linearized (semi)-implicit time discretizations.
\end{abstract}
\section{Introduction}
\setcounter{equation}{0}
Incompressible miscible flow in porous media is
governed in general by the following system of equations:
\begin{align}
&\Phi\frac{\partial c}{\partial t}-\nabla\cdot(D({\bf u})\nabla
c)+{\bf u}\cdot\nabla c= \hat c q^I-cq^P, \label{e-fuel-1}\\[3pt]
&\nabla\cdot\u=q^I-q^P, \label{e-fuel-2}\\[3pt]
&\u=-\frac{k(x)}{\mu(c)}\nabla p, \label{e-fuel-3}
\end{align}
where $p$ is the pressure of the fluid mixture, $\u$ is the velocity
and $c$ is the concentration of the first component; $k(x)$ is the
permeability of the medium, $\mu(c)$ is the concentration-dependent
viscosity, $\Phi$ is the porosity of the medium, $q^I$ and $q^P$ are
given injection and production sources, $\hat c$ is the
concentration of the first component in the injection source, and
$D({\bf u})=[D_{ij}({\bf u})]_{d\times d}$ is the
diffusion-dispersion tensor which may be given in different forms
(see \cite{BB1,BB2} for details). We assume that the system is
defined in a bounded smooth domain $\Omega$ in $\R^d$ $(d=2,3)$, for
$t\in[0,T]$, coupled with the initial and boundary conditions:
\begin{align}
\label{e-fuel-4}
\begin{array}{ll}
{\bf u}\cdot {\bf n}=0,~~
D({\bf u})\nabla c\cdot {\bf n}=0
&\mbox{for}~~x\in\partial\Omega,~~t\in[0,T],\\[3pt]
c(x,0)=c_0(x)~~ &\mbox{for}~~x\in\Omega.
\end{array}
\end{align}

The above system was investigated extensively in the last several
decades \cite{AE,CZW,ELSWZ, ERW, Feng} due to its wide
applications in various engineering areas, such as reservoir
simulations and exploration of underground water, oil and gas.
Existence of weak solutions of the system was obtained by Feng
\cite{Feng} for the 2D model and by Chen and Ewing \cite{ChenEwing}
for the 3D problem. Existence of semi-classical/classical solutions
is unknown so far. Numerical simulations have been done extensively
with various applications, see \cite{Dou,DFP,Peaceman} and the
references therein. Optimal error estimates of a standard
Galerkin-Galerkin method for the system
(\ref{e-fuel-1})-(\ref{e-fuel-4}) in two-dimensional space was
obtained first by Ewing and Wheeler \cite{EW} roughly under the
time-step condition $\tau =o(h)$. In their method, a linearized
semi-implicit Euler scheme was used in the time direction and
Galerkin FEM approximation was used both for the concentration and
the pressure. Later, a Galerkin-mixed finite element method was
proposed by Douglas et al \cite{DEW} for this system, where a
Galerkin approximation was applied for the concentration equation
and a mixed approximation in the Raviart--Thomas finite element
space was used for the pressure equation. A linearized
semi-implicit Euler scheme, the same as that used in \cite{EW},  was
applied for the time discretization. Optimal error estimates were
obtained under a similar time-step condition $\tau =o(h)$.
There are many other numerical
methods in literature for solving the equations of incompressible
miscible flow, such as \cite{Wang,WW} with an ELLAM in
two-dimensional space, \cite{KWang} with an MMOC-MFEM approximation
for the 2D problem, \cite{Duran,SY} with a characteristic-mixed
method in two and three dimensional spaces, respectively, and
\cite{MaNing,MLY} with a collocation-mixed method and a
characteristic-collocation method, respectively. In all those works,
error estimates were established under certain time-step conditions.
Moreover, linearized semi-implicit schemes have also been analyzed
with certain time-step restrictions for many other nonlinear parabolic-type
systems,
such as Navier-Stokes equations \cite{AG,He,KL,Liu1,Liu2},
nonlinear thermistor problems \cite{EL, Zhao},
viscoelastic fluid flow \cite{CLin,EM,WHS},
KdV equations \cite{MS},
nonlinear Schr\"{o}dinger equation \cite{S,Tou} and
some other equations \cite{ACM,HLS,SS}.
A key issue in analysis of FEMs is the boundedness of the
numerical solution in $L^{\infty}$ norm or a
stronger norm, which in a routine way can be estimated by
the mathematical induction with an inverse inequality,
such as,
\begin{equation}
\| u_h^n - R_hu(\cdot, t_n) \|_{L^\infty}
\le C h^{-d/2} \| u_h^n - R_h u(\cdot, t_n) \|_{L^2}
\le C h^{-d/2} (\tau^{m} + h^{r+1}) ,
\end{equation}
where $u^n_h$ is the finite element solution, $u$
is the exact solution, and $R_h$ is certain projection operator.
A time-step restriction arises immediately from the above inequality.
Such a time-step restriction may result in the use of a very small time step
and
extremely time-consuming in practical computations. The problem
becomes more serious when a non-uniform mesh is used.  However, we
believe that those time-step restrictions may not be necessary in most
cases.

In this paper, we analyze the linearized semi-implicit Euler scheme
with a popular Galerkin-mixed finite element approximation in the
spatial direction for the system (\ref{e-fuel-1})-(\ref{e-fuel-4}).
We establish optimal $L^2$ error estimates almost without any
time-step restriction (or when $h$ and $\tau$ are smaller than some
positive constants, respectively). Our theoretical
analysis is based on an error splitting proposed in \cite{LS} (also see
\cite{LiPhDThesis}) for a Joule heating system with a standard Galerkin FEM.
By introducing a corresponding time-discrete system, we split the numerical
error into two parts, the error in the temporal direction and the error in
the spatial direction.
Thus with the solution $v^n$ of the time-discrete system,
the numerical solution can be bounded by
\begin{equation}
\| u_h^n - R_h v^n \|_{L^\infty} \le C h^{-d/2} \| u_h^n - R_h v^n \|_{L^2}
\le C h^{-d/2} h^{r+1}
\end{equation}
without any time-step condition.
Rigorous analysis of the regularity of
the solution to the time-discrete PDEs is a key to our approach. With such proved
regularity, we obtain optimal and $\tau$-independent $L^2$-error
estimates of the Galerkin-mixed FEM for the time-discrete PDEs.
Our analysis presented in this paper provides a new understanding
on the commonly-used linearized schemes and clears up the misgivings
for the time-step size in practical computations.

The rest of the paper is organized as follows. In Section 2, we
introduce the linearized semi-implicit Euler scheme with a
Galerkin-mixed approximation in the spatial direction for the system
(\ref{e-fuel-1})-(\ref{e-fuel-4}) and present our main results. In
Section 3, we present a priori estimates of the solution to the
corresponding time-discrete system and the error estimate of the
linearized scheme. Optimal error estimates of the fully discrete
scheme in the $L^2$ norm are given in Section 4. In Section 5, we present
numerical examples to illustrate the convergence rate and the unconditional
convergence (stability) of the numerical method.
Conclusions are drawn in the last section.

\section{The Galerkin-mixed FEM and the main results}
\setcounter{equation}{0}
For any integer $m\geq 0$ and $1\leq p\leq\infty$, let $W^{m,p}$ be the
Sobolev space of functions defined on $\Omega$ equipped with the
norm
$$
\|f\|_{W^{m,p}}=\biggl(\sum_{|\beta|\leq m}\int_\Omega|D^\beta f|^p\d
x\biggl)^\frac{1}{p},
$$
where
$$
D^\beta=\frac{\partial^{|\beta|}}{\partial
x_1^{\beta_1}\cdots\partial x_d^{\beta_d}}
$$
for the multi-index $\beta=(\beta_1,\cdots,\beta_d)$, $\beta_1\geq
0$, $\cdots$, $\beta_d\geq 0$, and $|\beta|=\beta_1+\cdots+\beta_d$.
For any integer $m\geq 0$ and $0<\alpha<1$, let $C^{m+\alpha}$
denote the usual H\"{o}lder space equipped with the norm
$$
\|f\|_{C^{m+\alpha}}=\sum_{|\beta|\leq m}\|D^\beta
f\|_{C(\overline\Omega)}+\sum_{|\beta|=
k}\sup_{x,y\in\Omega}\frac{|D^\beta f(x)-D^\beta
f(y)|}{|x-y|^\alpha}.
$$
Let $I=(0,T)$. For any Banach space $X$, we consider functions
$g:I\rightarrow X$ and define the norm
$$
\|g\|_{L^p(I;X)} =\left\{
\begin{array}{ll}
\displaystyle\biggl(\int_0^T\|g(t)\|_X^pdt\biggl)^\frac{1}{p}, &
1\leq p<\infty
,\\[10pt]
\displaystyle{\rm ess\,sup}_{t\in I}\|g(t)\|_X, & p=\infty.
\end{array}
\right.
$$
In addition, we define $L^p_0$ as the subspace of $L^p$ consisting
of functions in $L^p$ whose integral over $\Omega$ are zeros.
Finally, we denote by $H$ the space of vector-valued functions
${\vec f}\in L^2\times L^2\times L^2$ such that $\nabla\cdot {\vec
f}\in L^2$.

Let $\pi_h$ be a regular division of $\Omega$ into triangles
$T_j$, $j=1,\cdots,M$, in $\R^2$
(or tetrahedra in $\R^3$), with $\Omega_h=
\cup_{j} T_j$ and denote by $h=\max_{1\leq j\leq
M}\{\mbox{diam}\,T_j\}$ the mesh size. For a triangle
$T_j$ at the boundary, we define $\widetilde T_j$ as the extension of $T_j$
to the triangle with one curved edge (or a tetrahedra with one curved face
in $\R^3$). For a given division of $\Omega$, we define the finite
element spaces (more details are described in \cite{LiPhDThesis}):
\begin{align*}
&S_h =\{w_h\in L^2(\Omega): w_h|_{\widetilde T_j}
\mbox{~is~linear~~for~each~element~$T_j\in \pi_h$ and~$\int_\Omega w_h
dx=0$}\},\\
&V_h =\{w_h\in C^0(\overline\Omega): w_h|_{\widetilde T_j}
\mbox{~is~linear~for~each~element~}  T_j\in \pi_h\} .
\end{align*}
Let $H_h$ be the subspace of $H$, as introduced by Raviart and Thomas
\cite{RT,VThomee} such that $\psi\cdot{\bf n}=0$ on $\partial\Omega$ and
${\rm div}\,\psi\in S_h$  for $\psi\in
H_h$.

In the rest part of this paper, we assume that the solution to the
initial-boundary value problem (\ref{e-fuel-1})-(\ref{e-fuel-4}) exists
and satisfies
\begin{align}
\label{StrongSOlEST}
&\|p\|_{L^\infty(I;H^3)}+\|{\bf u}\|_{L^\infty(I;H^2)}
+\|{\bf u}_t\|_{L^2(I;W^{1,3/2})}+\|c\|_{L^\infty(I;W^{2,s})} \nn \\
&+\|c_t\|_{L^\infty(I;H^2)}+\|c_t\|_{L^s(I;W^{1,s})}
+\|c_{tt}\|_{L^s(I;L^s)} \leq C
\end{align}
for some $s>d$ and
\begin{align}
\| q^I \|_{H^1}, \, \| q^P \|_{H^1} \leq C  .\label{qq}
\end{align}
Correspondingly, we assume that the permeability $k\in
W^{2,\infty}(\Omega)$ and satisfies
\begin{align}
k_0^{-1}\leq k(x)\leq k_0\quad\mbox{for}~~x\in\Omega,
\end{align}
the concentration-dependent viscosity $\mu\in C^1(\R)$ and
satisfies
\begin{align}
\mu_0^{-1}\leq \mu(s)\leq \mu_0\quad\mbox{for}~~s\in\R,
\end{align}
for some positive constant $\mu_0$. Moreover, we assume that the
diffusion-dispersion tensor is given by
$D({\bf u}) = \Phi  d_m I + D^*({\bf u})$,
where $d_m>0$, $D^*({\bf
u}) = d_1({\bf u}) I + d_2({\bf u})({\bf u} \otimes {\bf u})$ is
symmetric and positive definite and $\partial_{u_i} D\in
L^{\infty}$, $\partial_{u_i u_j}^2 D\in
L^{\infty}$ \cite{BB2}. For the initial-boundary value problem
(\ref{e-fuel-1})-(\ref{e-fuel-4}) to be well-posed, we require
\begin{equation}
\int_\Omega q^I \d x=\int_\Omega q^P \d x .
\end{equation}

Let $\{ t_n \}_{n=0}^N$ be a uniform partition of the time interval
$[0,T]$ with $\tau=T/N$ and denote
$$
p^n(x) = p(x,t_n),\quad {\bf u}^n(x)  = {\bf u}(x,t_n), \quad c^n(x)  = c(x,t_n) \, .
$$
For any sequence of functions $\{ f^n \}_{n=0}^N$, we define
$$
D_t f^{n+1}=\frac{f^{n+1}-f^n}{\tau}\, .
$$

The fully discrete mixed finite element scheme is to find $P_h^n\in
S_h/\{\rm constant\}$, $U_h^n\in H_h$ and ${\cal C}_h^n\in V_h$,
$n=0,1,\cdots,N$, such that for all $v_h\in H_h$, $\varphi_h\in S_h$
and $\phi_h\in V_h$,
\begin{align}
& \biggl(\frac{\mu({\cal C}^n_h)}{k(x)} U^{n+1}_h,\,v_h\biggl)
=-\Big(P^{n+1}_h ,\, \nabla \cdot v_h \Big),
\label{e-FEM-1}\\[3pt]
& \Big(\nabla\cdot U^{n+1}_h ,\, \varphi_h\Big)
=\Big(q^I-q^P,\, \varphi_h\Big),
\label{e-FEM-2}\\[3pt]
& \Big(\Phi  D_t{\cal C}^{n+1}_h, \, \phi_h\Big) + \Big(D(
U^{n+1}_h)
\nabla {\cal C}^{n+1}_h, \, \nabla \phi_h \Big) \nn\\
&~~~~~ + \Big( U^{n+1}_h\cdot\nabla {\cal C}^{n+1}_h,\,
\phi_h\Big)= \Big(\hat c q^I-{\cal C}^{n+1}_h q^P, \, \phi_h\Big)
\label{e-FEM-3}
\end{align}
where the initial data ${\cal C}_h^0$ is the Lagrangian piecewise
linear interpolation of $c^0$.

In this paper, we denote by $C$ a generic positive constant and
by $\epsilon$ a generic small positive constant, which
are independent of $n$, $h$ and $\tau$.
We present our main results in the following theorem.
\medskip
\begin{theorem}\label{ErrestFEMSol}
{\it Suppose that the initial-boundary value problem
{\rm (\ref{e-fuel-1})-(\ref{e-fuel-4})} has a unique solution $(p, {\bf u},c)$
which satisfies {\rm (\ref{StrongSOlEST})}. Then there exist positive constants
$h_0$ and $\tau_0$ such that when $h<h_0$ and $\tau<\tau_0$, the
finite element system {\rm (\ref{e-FEM-1})-(\ref{e-FEM-3})} admits a
unique solution $(P^n_h, U^n_h, {\cal C}^n_h)$, $n=1,\cdots,N$, which
satisfies
that
\begin{align}
&\max_{1\leq n\leq N}\|P^n_h - p^n\|_{L^2} +\max_{1\leq n\leq N}\|
U^n_h - {\bf u}^n\|_{L^2} +\max_{1\leq n\leq N}\|{\cal C}^n_h -c^n\|_{L^2}
\leq C(\tau+h^2).
\end{align}
}
\end{theorem}
\bigskip
We will present the proof of Theorem 2.1 in the next two sections.
The key to our proof is the following error splitting
\begin{align*}
& \|U^n_h- {\bf u}^n\|_{L^2} \le \| e^n_u \|_{L^2} + \| U^n - U^n_h
\|_{L^2},
\\
& \|P_h^n- p^n\|_{L^2} \le \| e^n_p \|_{L^2} + \| P^n - P_h^n
\|_{L^2},
\\
& \|{\cal C}^n_h -c^n \|_{L^2}  \le \| e^n_c \|_{L^2} + \| {\cal
C}^n - {\cal C}^n_h \|_{L^2},
\end{align*}
where
\begin{align*}
& e^n_p = P^n - p^n,  \\
& e^n_u = U^n - {\bf u}^n,  \\
& e^n_c = {\cal C}^n - c^n, \nn
\end{align*}
and $(P^n, U^n, {\cal C}^n)$ is the solution of a time-discrete
system defined in next section.
\section{Error estimates for time-discrete system}\label{SEction3}
\setcounter{equation}{0}
We define the time-discrete solution $(P^n, U^n, {\cal C}^n)$ by the
following elliptic system:
\begin{align}
&U^{n+1}=-\frac{k(x)}{\mu({\cal C}^n)}\nabla P^{n+1},
\label{TDe-fuel-1}
\\[3pt]
&\nabla\cdot U^{n+1}=q^I-q^P, \label{TDe-fuel-2}\\[3pt]
&\Phi D_t{\cal C}^{n+1} -\nabla\cdot(D( U^{n+1})\nabla {\cal
C}^{n+1}) + U^{n+1} \cdot\nabla{\cal C}^{n+1}= \hat c q^I-{\cal
C}^{n+1}q^P, \label{TDe-fuel-3}
\end{align}
for $x\in\Omega$ and $t\in[0,T]$, with the initial and boundary
conditions
\begin{align}
\label{TDBC}
\begin{array}{ll}
U^{n+1}\cdot {\bf n}=0,~~
D({U}^{n+1})\nabla {\cal C}^{n+1}\cdot {\bf n}=0
&\mbox{for}~~x\in\partial\Omega,~~t\in[0,T],\\[3pt]
{\cal C}^0(x)=c_0(x)~~ &\mbox{for}~~x\in\Omega,
\end{array}
\end{align}
The condition $\int_\Omega P^{n+1}dx =0 $ is enforced for the uniqueness of
solution.

In this section, we prove the existence, uniqueness and
regularity of the solution of the above
time-discrete system.
\begin{theorem}\label{ErrestTDSol}
{\it Suppose that the initial-boundary value problem
{\rm (\ref{e-fuel-1})-(\ref{e-fuel-4})} has a unique solution $(p, {\bf u},c)$
which satisfies {\rm (\ref{StrongSOlEST})}. Then there exists a positive constant
$\tau_0$ such that when $\tau<\tau_0$, the time-discrete system
{\rm (\ref{TDe-fuel-1})-(\ref{TDBC})} admits a unique solution $(P^n, U^n
, {\cal C}^n )$, $n=1,\cdots,N$, which satisfies
\begin{align}
\label{StrongSOlESTTD}
&\|P^{n}\|_{H^2} +\| U^n\|_{H^2}+\|{\cal
C}^n\|_{W^{2,s}}+\|D_t{\cal C}^n\|_{L^s}  +\|\nabla {\cal
C}^n\|_{L^\infty} \\
&+\biggl(\sum_{n=1}^{N}\tau\|
D_tU^{n}\|_{L^3}^2\biggl)^{\frac{1}{2}}
+\biggl(\sum_{n=1}^{N}\tau\|\nabla
D_tU^{n}\|_{L^{3/2}}^2\biggl)^{\frac{1}{2}}
+\biggl(\sum_{n=1}^{N}\tau\|
D_tC^{n}\|_{H^2}^2\biggl)^{\frac{1}{2}}
\leq C , \nn
\end{align}
and
\begin{align}\label{errorESTTD}
&\max_{1\leq n\leq N}\| e^n_p \|_{L^s} +\max_{1\leq n\leq N}\| e_u^n
\|_{L^s} +\max_{1\leq n\leq N}\| e_c^n \|_{L^s}
\leq C \tau .
\end{align}
}
\end{theorem}
\medskip
\noindent{\it Proof}~~~ It suffices to establish the estimates
presented in (\ref{StrongSOlESTTD}). With such estimates, existence
and uniqueness of solution follow a routine way.
We observe that $e^n_p$, $e_u^n$ and $e_c^n$ satisfy the following
equations
\begin{align}
&-\nabla\cdot\biggl(\frac{k(x)}{\mu({\cal C}^n)}\nabla e^{n+1}_p\biggl)
= \nabla\cdot\biggl[\biggl(\frac{k(x)}{\mu({\cal C}^{n})}
-\frac{k(x)}{\mu(c^{n})}\biggl)\nabla p^{n+1}\biggl],
\label{TDerr-fuel-1}
\\[10pt]
& e_u^{n+1}=-\frac{k(x)}{\mu({\cal C}^{n})}\nabla e^{n+1}_p  -
\biggl(\frac{k(x)}{\mu({\cal C}^{n})}
-\frac{k(x)}{\mu(c^{n})}\biggl)\nabla p^{n+1},
\label{TDerr-fuel-2}
\\[10pt]
&\Phi D_t e_c^{n+1} -\nabla\cdot(D(U^{n+1})\nabla e_c^{n+1}) +
U^{n+1}\cdot\nabla e_c^{n+1}
\nn\\
&=\nabla\cdot\Big((D( U^{n+1})-D({\bf u}^{n+1}))\nabla c^{n+1}\Big)
-(U^{n+1}-{\bf u}^{n+1})\cdot\nabla c^{n+1}-e_c^{n+1}q^P +{\cal E}^{n+1},
\label{TDerr-fuel-3}
\end{align}
for $x\in\Omega$ and $t\in[0,T]$, with the initial and boundary
conditions
\begin{align}
\label{D-BC}
\begin{array}{ll}
\displaystyle\frac{k(x)}{\mu({\cal C}^{n})} \nabla e_p^{n+1}\cdot
{\bf n}=0,~~ D(U^{n+1})\nabla e_c^{n+1}\cdot {\bf n}=0,
&~~\mbox{for}~~x\in\partial\Omega,~~t\in[0,T],\\[10pt]
e_c^0(x)=0,~~ &~~\mbox{for}~~x\in\Omega \, ,
\end{array}
\end{align} where
$$
{\cal E}^{n+1} = \Phi (c_t^{n+1}-D_t c^{n+1})
$$
is the truncation error due to the discretization in the time direction.
From the regularity assumption for $c$ in (\ref{StrongSOlEST}), we
see that
\begin{align}\label{TDtruncerr}
&\|{\cal E}^{n+1}\|_{L^2}\leq C,\quad \sum_{n=0}^{N-1}\tau \|{\cal
E}^{n+1}\|_{L^2}^2 \leq C\tau^2  ,\\
&\|{\cal E}^{n+1}\|_{L^s}\leq C,\quad \sum_{n=0}^{N-1}\tau \|{\cal
E}^{n+1}\|_{L^s}^s \leq C\tau^s .
\label{TDtruncerr123}
\end{align}

To prove the error estimate (\ref{errorESTTD}), we multiply
(\ref{TDerr-fuel-1}) by $e_p^{n+1}$ to get
\begin{align}
\|\nabla e_p^{n+1}\|_{L^2} &\leq \biggl
\|\biggl(\frac{k(x)}{\mu({\cal C}^{n})}
-\frac{k(x)}{\mu(c^{n})}\biggl)\nabla p^{n+1} \biggl \|_{L^2} \le C
\|e_c^{n}\|_{L^2} \| \nabla p^{n+1} \|_{L^\infty} \leq C
\|e_c^{n}\|_{L^2} \label{TDerr1fuel-1}
\end{align}
which together with (\ref{TDerr-fuel-2}) implies that
\begin{align}
\|e_u^{n+1}\|_{L^2} & \leq C\|\nabla e_p^{n+1}\|_{L^2} +
C\|e_c^{n}\|_{L^2} \| \nabla p^{n+1} \|_{L^\infty} \le  C
 \|e_c^{n}\|_{L^2} \, . \label{TDerr1fuel-2}
\end{align}
Then we multiply (\ref{TDerr-fuel-3}) by $e_c^{n+1}$ to get
\begin{align}
\frac{1}{2}& D_t \biggl( \Phi \|e_c^{n+1}\|_{L^2}^2 \biggl)
 + \| \sqrt{D(U^{n+1})}\nabla e_c^{n+1} \|_{L^2}^2
\nn \\
& \leq C\| e^{n+1} \|_{L^4}^2 \|q^I-q^P\|_{L^2}+C\|e_u^{n+1}\|_{L^2} \|
\nabla e_c^{n+1} \|_{L^2} \|\nabla c^{n+1}\|_{L^\infty}
\nn \\
&~~ + C\|e_u^{n+1}\|_{L^2} \|e_c^{n+1}\|_{L^2} \|\nabla
c^{n+1}\|_{L^\infty}+ C \|e_c^{n+1}\|_{L^4}^2 \| q^P \|_{L^2}  +
C\|{\cal E}^{n+1}\|_{L^2} \| e_c^{n+1} \|_{H^1}
\nn \\
& \leq C\|e_c^{n+1}\|_{L^4}^2 + C \| e_c^{n+1} \|_{H^1}
\|e_c^{n}\|_{L^2}  + C\epsilon^{-1}\|{\cal
E}^{n+1}\|_{L^2}^2+\epsilon \| e_c^{n+1} \|_{H^1}^2
\nn \\
& \leq \epsilon  \| \nabla e_c^{n+1}  \|_{L^2}^2 + C_\epsilon
\|e_c^{n+1}\|_{L^2}^2 + C \|e_c^n\|_{L^2}^2 + C\|{\cal
E}^{n+1}\|_{L^2}^2 \label{TDerr1-fuel-3} ,
\end{align}
where we have used the inequalities
$$
| ( U^{n+1} \cdot \nabla e_c^{n+1}, \,  e_c^{n+1} )| =  | ( \nabla
\cdot U^{n+1}, \, |e_c^{n+1}|^2 )| \leq  \| e^{n+1} \|_{L^4}^2
\|q^I-q^P\|_{L^2},
$$
and
$$
\| e_c^{n+1} \|_{L^4}^2 \leq \epsilon \| \nabla e_c^{n+1} \|_{L^2}^2
+ C_\epsilon \| e_c^{n+1} \|_{L^2}^2 \, .
$$
The square root of $D(U^{n+1})$ exists because $D(U^{n+1})$ is a
symmetric and positive definite matrix. With Gronwall's inequality
and (\ref{TDtruncerr}), (\ref{TDerr1-fuel-3}) reduces to
\begin{align}\label{fdsmkwjerqio984}
\|e_c^{n+1}\|_{L^2}^2 + \sum_{n=0}^{N-1}\tau \big \|\nabla e_c^{n+1}
\big \|_{L^2}^2 \leq C\tau^2 ,
\end{align}
provided $\tau<\tau_1$ for some positive constant $\tau_1$. The above
inequality implies that
\begin{align}\label{fdnhifho08}
\| D_t e_c^{n+1} \|_{L^2} \le C
\, .
\end{align}
From (\ref{TDerr1fuel-1})-(\ref{TDerr1fuel-2}), we derive that
\begin{align}\label{safnuyoe90}
\| e_u^{n+1} \|_{L^2} + \| \nabla e_p^{n+1} \|_{L^2} \le C \tau \, .
\end{align}
To prove (\ref{StrongSOlESTTD}), we rewrite (\ref{TDerr-fuel-3}) as
\begin{align}
&\Phi D_t e_c^{n+1} -\nabla\cdot(D( {\bf u}^{n+1})\nabla e_c^{n+1})
\\
&=-\nabla\cdot\Big((D( {\bf u}^{n+1})-D(U^{n+1}))\nabla
e_c^{n+1}\Big)-(U^{n+1}-{\bf u}^{n+1})\cdot\nabla e_c^{n+1}-{\bf u}^{n+1}\cdot\nabla
e_c^{n+1}
\nonumber\\
&~~~~+\nabla\cdot\Big((D( U^{n+1})-D( {\bf u}^{n+1}))\nabla
c^{n+1}\Big)-(U^{n+1}-{\bf u}^{n+1})\cdot\nabla c^{n+1}-e_c^{n+1}q^P +{\cal E}^{n+1}
\, .\nonumber
\end{align}
Since ${\bf u}^{n+1} \in C^{\alpha}(\overline\Omega)$, using the $W^{1,6}$
estimates of elliptic equations for $e_c$ \cite{ChenYZ}, we get
\begin{align*}
\|\nabla e_c^{n+1}\|_{L^6} & \leq C \| D_t e_c^{n+1} \|_{L^2} + C
\|e_u^{n+1} \|_{L^\infty}\|\nabla
e_c^{n+1}\|_{L^6}+C\|{\bf u}^{n+1}\|_{L^\infty}\|\nabla e_c^{n+1}\|_{L^2}\\
&~~ + C\|e_u^{n+1} \|_{L^\infty}\|\nabla c^{n+1}\|_{L^6} +
C\|q^P\|_{L^6} \| e_c^{n+1} \|_{L^3} + C\| {\cal E}^{n+1} \|_{L^2}
\end{align*}
which further reduces to
\begin{align}\label{errestecw16}
\| \nabla e_c^{n+1}\|_{L^6} \leq  C_0\|e_u^{n+1} \|_{L^\infty} ( C_0
+ \|\nabla e_c^{n+1}\|_{L^6}) + C_0
\end{align}
for some positive constant $C_0$ independent of $n,\tau,h$.

On the other hand,  we can rewrite the equation (\ref{TDerr-fuel-1}) into
the following form:
\begin{align}
&-\frac{k(x)}{\mu({\cal C}^n)}\Delta e_p^{n+1} =
\nabla\biggl(\frac{k(x)}{\mu({\cal C}^n)}\biggl) \cdot\nabla
e_p^{n+1} + \nabla\cdot\biggl[\biggl(\frac{k(x)}{\mu({\cal C}^n)}
-\frac{k(x)}{\mu(c^n)}\biggl)\nabla p^{n+1}\biggl] . \label{migenwo}
\end{align}
Since we have
\begin{align*}
\biggl\|\nabla\biggl(\frac{k(x)}{\mu({\cal C}^n)}\biggl) \cdot\nabla
e_p^{n+1}\biggl\|_{L^6} & \leq
\biggl\|\nabla\biggl(\frac{k(x)}{\mu({\cal C}^n)}\biggl) \biggl
\|_{L^6} \| \nabla e_p^{n+1} \|_{L^\infty}
\\
& \leq (C+C\|\nabla{\cal C}^n\|_{L^6})\|\nabla
e_p^{n+1}\|_{L^{\infty}}
\end{align*}
and
\begin{align*}
& \biggl \|
 \nabla\cdot\biggl[ \biggl(\frac{k(x)}{\mu({\cal C}^n)}
-\frac{k(x)}{\mu(c^n)}\biggl) \nabla p^{n+1}\biggl] \biggl \|_{L^6}
\\
& \le \biggl\| \nabla\biggl( \frac{k(x)}{\mu({\cal C}^n)}
-\frac{k(x)}{\mu(c^n)}\biggl) \cdot \nabla p^{n+1} \biggl \|_{L^6} +
\biggl\|\biggl( \frac{k(x)}{\mu({\cal C}^n)}
-\frac{k(x)}{\mu(c^n)}\biggl) \nabla \cdot \nabla p^{n+1} \biggl
\|_{L^6}
\\
& \le \| \nabla e_c^n \|_{L^6} \| \nabla p^{n+1} \|_{L^{\infty}} +
\| e_c^n \|_{L^{\infty}} \| p^{n+1} \|_{W^{2,6}} ,
\end{align*}
by the $W^{2,6}$ estimates of elliptic equations \cite{ChenYZ}, we
derive from (\ref{migenwo}) that
\begin{align}
\| e_p^{n+1} \|_{W^{2,6}} & \leq (C+C\|\nabla {\cal
C}^n\|_{L^6})\|\nabla e_p^{n+1}\|_{L^\infty} + C \| \nabla e_c^n
\|_{L^6} + C \| e_c^n \|_{L^\infty}
\nn \\
& \leq (C+C\|\nabla e_c^n\|_{L^6}) ( \epsilon \|
e_p^{n+1}\|_{W^{2,6}} + C_{\epsilon}\| e_p \|_{L^2}) \nn \\
&~~~ + (C+\epsilon) \| \nabla e_c^n \|_{L^6} + (C+C_\epsilon)\|
e_c^n \|_{L^2}
\nn \\
& \leq (C+ C\|\nabla e_c^n\|_{L^6} )\epsilon \|
e_p^{n+1}\|_{W^{2,6}} + C_{\epsilon}+ C_{\epsilon}\|\nabla
e_c^n\|_{L^6}  \label{p-w26}
\end{align}
From (\ref{TDerr-fuel-2}) we observe that
\begin{align}
\|e_u^{n+1}\|_{W^{1,6}}
&
\leq (C+C\|e^{n}_c\|_{W^{1,6}})(\|\nabla
e^{n+1}_p\|_{L^\infty}+\|\nabla
p^{n+1}\|_{L^\infty})+C\|e^{n+1}_p\|_{W^{2,6}} \nn\\
& \leq (C+ C\|\nabla e_c^n\|_{L^6} )\epsilon \|
e_p^{n+1}\|_{W^{2,6}} + C_{\epsilon}+ C_{\epsilon}\|\nabla
e_c^n\|_{L^6} , \label{p-w26667}
\end{align}
and by the Sobolev interpolation inequality, we have
\begin{align}
\|e^{n+1}_u\|_{L^\infty}\leq
C\|e^{n+1}_u\|_{L^2}^{1/4}\|e^{n+1}_u\|_{W^{1,6}}^{3/4} \leq
C\tau^{1/4}\|e^{n+1}_u\|_{W^{1,6}}^{3/4}. \label{p-w2666}
\end{align}

With the estimates (\ref{errestecw16}), (\ref{p-w26}),
(\ref{p-w26667}) and (\ref{p-w2666}), we now apply mathematical
induction to prove
\begin{align}\label{nkldafihylewat}
\|\nabla e^n_c\|_{L^6}\leq 4C_0
\end{align}
for $n=1,2,\cdots,N$. Clearly, the above inequality holds when
$n=0$. With the above inequalities, by choosing a proper small
$\epsilon=\epsilon(C_0)$, we derive from
(\ref{p-w26})-(\ref{p-w2666}) that
\begin{align}
\|e^{n+1}_u\|_{L^\infty}\leq C_1\tau^{1/4}, \label{p-w266666}
\end{align}
where $C_1$ may depend on $C_0$. Hence, there exists a positive
constant $\tau_0=\tau_0(C_1)$ such that when $\tau<\tau_0$, we have
$$
C_0\|e^{n+1}_u\|_{L^\infty}\leq 1/2 .
$$
Substituting the above inequality into (\ref{errestecw16}) gives
$$
\|\nabla e^{n+1}_c\|_{L^6}\leq 4C_0 .
$$
By mathematical induction, (\ref{nkldafihylewat}) holds for
$1\leq n\leq N$. From (\ref{p-w26}) and (\ref{p-w26667}), we also
get
\begin{align}
&\max_{1\leq n\leq N}\| e_p^{n+1} \|_{W^{2,6}}+\max_{1\leq n\leq
N}\| e_u^{n+1} \|_{W^{1,6}} \leq C . \label{pW26bound}
\end{align}

By a similar approach as (\ref{p-w26667}), we can prove that $\|\nabla
e^{n+1}_u\|_{L^{3/2}}\leq C(\| e^{n}_c\|_{L^2}
+\|\nabla e^{n}_c\|_{L^{3/2}})$ and so
\begin{align*}
&\|\nabla D_tU^{n+1}\|_{L^{3/2}}\leq \|\nabla
D_te_u^{n+1}\|_{L^{3/2}}+\|\nabla D_t{\bf u}^{n+1}\|_{L^{3/2}}\\
&\leq C\tau^{-1}\|\nabla e_u^{n+1}\|_{L^{3/2}}+C\tau^{-1}\|\nabla
e^n_u\|_{L^{3/2}}+\|\nabla D_t{\bf u}^{n+1}\|_{L^{3/2}}\\
&\leq C\tau^{-1}(\| e^{n}_c\|_{L^2}
+\|\nabla e^{n}_c\|_{L^{3/2}})+C\tau^{-1}(\| e^{n-1}_c\|_{L^2}
+\|\nabla e^{n-1}_c\|_{L^{3/2}})+\|\nabla D_t{\bf u}^{n+1}\|_{L^{3/2}},
\end{align*}
and from (\ref{fdsmkwjerqio984}) we see that
\begin{align}\label{fdsjporop90705}
\biggl(\sum_{n=1}^{N-1}\tau\|\nabla
D_tU^{n+1}\|_{L^{3/2}}^2\biggl)^{\frac{1}{2}}\leq C .
\end{align}
Since $W^{1,3/2}\hookrightarrow L^3$ in $\R^d$ ($d=2,3$), it follows that
$$\|D_tU^{n+1}\|_{L^3}\leq C\|D_tU^{n+1}\|_{L^2}+C\|\nabla
D_tU^{n+1}\|_{L^{3/2}}.
$$
From (\ref{safnuyoe90}) and (\ref{fdsjporop90705}), we see that
\begin{align}\label{fdsjporop90704}
\biggl(\sum_{n=1}^{N-1}\tau\|
D_tU^{n+1}\|_{L^3}^2\biggl)^{\frac{1}{2}}\leq C .
\end{align}

For the $H^2$ regularity of ${\cal C}^{n+1}$, we rewrite
(\ref{TDe-fuel-3}) as
\begin{align}\label{dnidui88}
-\nabla\cdot(D( U^{n+1})\nabla {\cal C}^{n+1}) + U^{n+1}
\cdot\nabla{\cal C}^{n+1}= -\Phi D_t{\cal C}^{n+1} +\hat c q^I-{\cal
C}^{n+1}q^P.
\end{align}
With (\ref{fdnhifho08}), (\ref{nkldafihylewat}) and (\ref{pW26bound}), we
see that the right-hand side above is bounded in $L^2$ and so we can apply
the $H^2$
estimates for the elliptic equation \cite{Evans} to obtain
\begin{align}
\|{\cal C}^{n+1}\|_{H^2}\leq C.
\end{align}

We rewrite the equations (\ref{TDe-fuel-1})-(\ref{TDe-fuel-2}) as
\begin{align*}
-\nabla\cdot\biggl(\frac{k(x)}{\mu(C^n)}\nabla P^{n+1}\biggl)=q^I-q^P
\end{align*}
and apply the
$H^3$ estimates of
elliptic equations \cite{Evans} to the above equation. Then we get
\begin{align}\label{pH3uH2bound}
\|e_p^{n+1}\|_{H^3}+\|e_u^{n+1}\|_{H^2} \leq C .
\end{align}

Note that $H^2\hookrightarrow C^\alpha$ for some $\alpha>0$. With
the H\"{o}lder regularity of ${\cal C}^n$, applying the $W^{1,s}$
estimates to (\ref{TDerr-fuel-1}) and using (\ref{TDerr-fuel-2}), it
is not difficult to see that
\begin{align}\label{Lsesteunecn}
\|e_u^{n+1}\|_{L^s}  \leq  C
 \|e_c^{n}\|_{L^s} \, .
\end{align}
Multiplying (\ref{TDerr-fuel-3}) by
$|e_c^{n+1}|^{s-2}e_c^{n+1}$ and using (\ref{Lsesteunecn}), one can
derive that
\begin{align}
&\int_\Omega s^{-1}\Phi D_t |e_c^{n+1}|^s \d x+
(s-1)\int_\Omega|e_c^{n+1}|^{s-2}D(U^{n+1})\nabla
e_c^{n+1}\cdot\nabla e_c^{n+1}\d x
\nn\\
&\leq\int_\Omega \int_\Omega s^{-1}(q^I-q^P)|e_c^{n+1} |^s\d x
+\|\nabla
c^{n+1}\|_{L^\infty}\|U^{n+1}-{\bf u}^{n+1}\|_{L^s}\|e_c^{n+1}\|_{L^s}^{s-1} \nn\\
&~~~+\|q^P\|_{L^\infty}\|e_c^{n+1}\|_{L^s}
+\|{\cal E}^{n+1}\|_{W^{-1,s}}\|e_c^{n+1}\|_{L^s}^{\frac{s-2}{2}}
\||e_c^{n+1}|^{s-2}\nabla e_c^{n+1}\|_{L^2}\nn \\
&~~~+(s-1)\|D( U^{n+1})-D({\bf u}^{n+1})\|_{L^s}\|\nabla
c^{n+1}\|_{L^\infty}\|e_c^{n+1}\|_{L^s}^{\frac{s-2}{2}}
\||e_c^{n+1}|^{s-2}\nabla e_c^{n+1}\|_{L^2} \nn\\
&\leq (C+\|{\cal
E}^{n+1}\|_{W^{-1,s}}^s)\epsilon^{-1}(\|e_c^{n+1}\|_{L^s}^s+\|e_c^n\|_{L^s}^s)+\epsilon
(s-1)\int_\Omega|e_c^{n+1}|^{s-2}|\nabla e_c^{n+1}|^2\d x . \nn
\end{align}
Choosing a proper small $\epsilon$ and using Gronwall's inequality
lead to
\begin{align*}
\|e_c^{n+1}\|_{L^s}  \leq C \tau ,
\end{align*}
where we have used (\ref{TDtruncerr123}).
It follows that
\begin{align}
\|D_te_c^{n+1}\|_{L^s}  \leq C  .
\end{align}
With the above estimate, by applying the $W^{2,s}$ estimate to
(\ref{dnidui88}), we obtain
\begin{align}
\|e_c^{n+1}\|_{W^{2,s}}  \leq C
\end{align}
and by the Sobolev embedding theorem we get $\|\nabla
e_c^{n+1}\|_{L^\infty}\leq C \|e_c^{n+1}\|_{W^{2,s}}  \leq C$.
By applying the $W^{1,s}$ estimates to the elliptic
equation (\ref{TDerr-fuel-1}) and using (\ref{TDerr-fuel-2}), we
obtain
\begin{align}
\|\nabla e_p^{n+1}\|_{L^s} +\|e_u^{n+1}\|_{L^s}  \leq C\tau .
\end{align}

Finally, we multiply (\ref{TDerr-fuel-3}) by $-\nabla\cdot(D(U^{n+1})\nabla
e^{n+1}_c)$ and summing up the results for $n=0,1,\cdots,N-1$. Then we get
\begin{align}
&D_t\big(D(U^{n+1})\nabla e^{n+1}_c,\nabla
e^{n+1}_c\big)+\big\|\nabla\cdot\big(D(U^{n+1})\nabla
e^{n+1}_c\big)\big\|_{L^2}^2 \label{chap10fsnioqo1236}\\
&\leq C\|D_tD(U^{n+1})\|_{L^3}\|\nabla e^{n+1}_c\|_{L^2}\|\nabla
e^{n+1}_c\|_{L^6}
+C\|e_u^{n+1}\|_{H^1}^2\|\nabla c^{n+1}\|_{L^\infty}^2 \nn\\
&~~~+C\|e_u^{n+1}\|_{L^6}^2\|c^{n+1}\|_{W^{2,3}}^2
+C\|e^{n+1}_c\|_{H^1}^2+C\|{\cal E}^{n+1}\|_{L^2}^2 \nn\\
&\leq C_\epsilon
\|D_tD(U^{n+1})\|_{L^3}^2\|\nabla e^{n+1}_c\|_{L^2}^2+\epsilon
\|e_c^{n+1}\|_{W^{1,6}}^2
+C(\|e_u^{n+1}\|_{H^1}^2+
\|e^{n+1}_c\|_{H^1}^2+\|{\cal E}^{n+1}\|_{L^2}^2) . \nn
\end{align}
Since by the Sobolev inequality and the theory of elliptic equations we have
\begin{align*}
&\|e^{n+1}_c\|_{W^{1,6}}\leq C\|e^{n+1}_c\|_{H^2}\leq
C\big\|\nabla\cdot\big(D(U^{n+1})\nabla e^{n+1}_c\big)\big\|_{L^2} ,
\end{align*}
and by the $H^1$ estimates of the equation (\ref{TDerr-fuel-1}) we have
\begin{align*}
&\|e^{n+1}_u\|_{H^1}\leq C\|e^{n+1}_c\|_{H^1}\|\nabla
p^{n+1}\|_{L^\infty}+C\|e^{n+1}_c\|_{L^3}\|p^{n+1}\|_{W^{2,6}}\leq
C\|e^{n+1}_c\|_{H^1} ,
\end{align*}
by choosing a small $\epsilon$, the inequality (\ref{chap10fsnioqo1236})
reduces to
\begin{align*}
&D_t\big(D(U^{n+1})\nabla e^{n+1}_c,\nabla e^{n+1}_c\big)+\frac{1}{2}\big
\|\nabla\cdot\big(D(U^{n+1})\nabla e^{n+1}_c\big)\big\|_{L^2}^2\\
&\leq C\|D_tD(U^{n+1})\|_{L^3}^2\|\nabla e^{n+1}_c\|_{L^2}^2
+C(\|\nabla e^{n+1}_c\|_{L^2}^2+\|e^{n+1}_c\|_{L^2}^2
+\|{\cal E}^{n+1}\|_{L^2}^2) \\
&\leq (C\|D_tD(U^{n+1})\|_{L^3}^2+C)
\big(D(U^{n+1})\nabla e^{n+1}_c,\nabla e^{n+1}_c\big)
+C(\|e^{n+1}_c\|_{L^2}^2+\|{\cal E}^{n+1}\|_{L^2}^2) .
\end{align*}
By applying Gronwall's inequality, using (\ref{TDtruncerr}),
(\ref{fdsmkwjerqio984}) and (\ref{fdsjporop90704}), we get
\begin{align*}
&\max_{1\leq n\leq N}\|\nabla
e^n_c\|_{L^2}^2+\sum_{n=1}^N\tau\|e^n_c\|_{H^2}^2 \leq C\tau^2 .
\end{align*}
In particular, the above inequality implies that
\begin{align*}
&\sum_{n=1}^N\tau\|D_t{\cal C}^n\|_{H^2}^2\leq\sum_{n=1}^N
\tau\|D_tc^n\|_{H^2}^2
+\sum_{n=1}^N\tau\|D_te^n_c\|_{H^2}^2 \leq
C+C\tau^{-2}\sum_{n=1}^N\tau\|e^n_c\|_{H^2}^2 \leq C .
\end{align*}

The proof of Theorem \ref{ErrestTDSol} is complete. \endproof
\section{Error estimates of the fully-discrete system}
\setcounter{equation}{0} To provide optimal error estimates for the fully
discrete scheme (\ref{e-FEM-1})-(\ref{e-FEM-3}), we define three
projections below.

Let $\Pi_h: L^2(\Omega)\rightarrow S_h$ be the $L^2$ projection defined by
$$
(\Pi_h \phi,\chi)=(\phi,\chi),~~\mbox{for all}~\phi\in
L^2~\mbox{and}~\chi\in S_h.
$$
For any fixed integer $n\geq 1$, let $\Pi_h^n: H^1(\Omega) \rightarrow
V_h$ be a projection defined by the following elliptic problem,
\begin{align}\label{njkleateh}
\Big(D(U^n)\nabla (v-\Pi_h^nv), \, \nabla \phi_h \Big) = 0 ,\quad
\mbox{for~all}~~\phi_h\in V_h,~~v\in H^1(\Omega)
\end{align}
with $\int_\Omega(v-\Pi_h^nv)d x=0$, and we define $\Pi_h^0:=\Pi_h^1$.
Moreover, let $Q_h: H \rightarrow H_h$ be a
projection such that \cite{VThomee}
\begin{align}
\Big(\nabla\cdot(w-Q_hw)\, , \chi \Big) = 0 ,\quad
\mbox{for~all}~~\chi\in S_h, ~ w\in H . \label{Q}
\end{align}
with a slight modification on the triangles/tetrahadons near the boundary.

By the theory of Galerkin and mixed finite element methods for
linear elliptic problems \cite{RS,VThomee}, with the
regularity $U^n\in H^2(\Omega)$, we have
\begin{align}
&\|v-\Pi_hv\|_{L^2}+h\|v-\Pi_hv\|_{H^1}\leq
Ch^2\|v\|_{H^2},\qquad\quad\,\mbox{for~all~}v\in H^2(\Omega),
\nn \\
&\|v-\Pi_h^{n+1}v\|_{L^p}+h\|v-\Pi_h^{n+1}v\|_{W^{1,p}}\leq
Ch^2\|v\|_{W^{2,p}},\quad \mbox{for~all~}v\in W^{2,p}(\Omega) ,
\nn\\
& \|w-Q_h w\|_{L^2}+h\|w-Q_h w\|_{H^1}\leq
Ch^2\|w\|_{H^2},\quad~\mbox{for~all~} w \in H,  \label{p-app}
\end{align}

Previous works on Galerkin (or mixed) FEM for the nonlinear
parabolic system (\ref{e-fuel-1})-(\ref{e-fuel-3}) required the
estimate
\begin{align}\label{knldfahlae0}
&\|\partial_t(c^{n+1} -\widetilde\Pi_h^{n+1}c^{n+1})\|_{L^2}\leq
Ch^2
\end{align}
for an elliptic projection operator $\widetilde\Pi_h^{n+1}$ defined
by
\begin{align}\label{njkleateh}
\Big(D({\bf u})\nabla (v-\widetilde\Pi_hv), \, \nabla \phi_h \Big) =
0 ,\quad \mbox{for~all}~~\phi_h\in V_h,~~v\in H^1(\Omega) .
\end{align}
The inequality (\ref{knldfahlae0}) was proved by Wheeler \cite{Whe}
based on the regularity assumption $\|\nabla D({\bf
u})_t\|_{L^\infty}\leq C$. In the following lemma, we will prove an
analogues inequality:
\begin{align}\label{knldfahlae}
&\biggl(\sum_{n=0}^{N-1}\tau\|D_t({\cal C}^{n+1}
-\Pi_h^{n+1}{\cal C}^{n+1})\|_{H^{-1}}^2\biggl)^{1/2}
\leq Ch^2,
\end{align}
based on weaker regularity of $U^n$ proved in the last section. The
above inequality is necessary to obtain optimal $L^2$ error
estimates. \vskip0.1in
\begin{lemma}\label{njekfqhklad}
{\it With the regularity assumption {\rm (\ref{StrongSOlEST})}, the
regularity of $({\cal C}^{n+1}, U^{n+1})$ given in
{\rm (\ref{StrongSOlESTTD})}, and the error estimates given in
{\rm (\ref{errorESTTD})}, the estimate {\rm (\ref{knldfahlae})} holds.
}
\end{lemma}
\noindent{\it Proof}~~ We only prove the lemma for the 3D problem
(with $s>3$ in Theorem \ref{ErrestTDSol}). The 2D problem can be
handled similarly. Note that
\begin{align}
&\Big(D(U^{n+1})\nabla ({\cal C}^{n+1}-\Pi_h^{n+1}{\cal C}^{n+1}),
\,
\nabla \phi_h \Big) = 0 ,\\
&\Big(D(U^n)\nabla ({\cal C}^{n+1}-\Pi_h^n{\cal C}^{n+1}), \, \nabla
\phi_h \Big) = 0 .
\end{align}
The difference of the above two equations gives
\begin{align*}
&\Big(D(U^{n+1})\nabla (\Pi_h^n{\cal C}^{n+1}-\Pi_h^{n+1}{\cal C}^{n+1}),
\, \nabla \phi_h \Big) \\
& + \Big ( (D(U^{n+1})-D(U^n))\nabla ({\cal C}^{n+1}-\Pi_h^n{\cal C}^{n+1}),
\, \nabla \phi_h \Big) =0 .
\end{align*}
By the $W^{1,p}$ estimates of elliptic projections \cite{RS}, we
have
\begin{align*}
&\|\nabla (\Pi_h^n{\cal C}^{n+1}-\Pi_h^{n+1}{\cal
C}^{n+1})\|_{L^{6/5}}\leq C \|(D(U^{n+1})-D(U^n))\nabla ({\cal
C}^{n+1}-\Pi_h^n{\cal C}^{n+1})\|_{L^{6/5}}\\
&\leq C \|D(U^{n+1})-D(U^n)\|_{L^2}\|\nabla ({\cal
C}^{n+1}-\Pi_h^n{\cal C}^{n+1})\|_{L^3} \leq
C\|D(U^{n+1})-D(U^n)\|_{L^2} h.
\end{align*}
For any $\varphi\in H^1(\Omega)$ with $\int_\Omega\varphi \d x=0$,
let $\psi$ be the solution of the equation
$$
-\nabla\cdot\Big(D(U^{n+1})\nabla \psi \Big)=\varphi
$$
with the boundary condition $\nabla\psi\cdot{\bf n}=0$ on
$\partial\Omega$. Since $U^{n+1}$ is uniformly bounded in
$H^2(\Omega)$, i.e. $\|U^{n+1}\|_{H^2}\leq C$, it is easy to check
that
$$
\|\psi\|_{H^3}\leq C\|\varphi\|_{H^1} .
$$
By the boundary condition $U^{n+1} \cdot{\bf n} = U^{n} \cdot{\bf
n}= 0$ and the expression of the function $D({\bf u})$, we see that
$D(U^{n+1})\nabla\psi\cdot{\bf n}=D(U^n)\nabla\psi\cdot{\bf n}=0$ on
$\partial\Omega$. We see that for $n\geq 1$,
\begin{align*}
&\big(\Pi_h^n{\cal C}^{n+1}-\Pi_h^{n+1}{\cal
C}^{n+1},\varphi)\\
&=\Big(D(U^{n+1})\nabla (\Pi_h^n{\cal C}^{n+1}-\Pi_h^{n+1}{\cal
C}^{n+1}), \, \nabla \psi \Big) \\
&=\Big(D(U^{n+1})
\nabla (\Pi_h^n{\cal C}^{n+1}-\Pi_h^{n+1}{\cal C}^{n+1}),
\, \nabla (\psi-P_h\psi) \Big)\\
&~~~ - \Big((D(U^{n+1})-D(U^n)) \nabla ({\cal C}^{n+1}-\Pi_h^n{\cal
C}^{n+1}),
\, \nabla (P_h\psi-\psi) \Big) \\
&~~~ - \Big( (D(U^{n+1})-D(U^n))
\nabla ({\cal C}^{n+1}-\Pi_h^n{\cal C}^{n+1}), \, \nabla \psi \Big) \\
&\leq C\|D(U^{n+1})-D(U^n)\|_{L^2} h\|\nabla
(\psi-P_h\psi)\|_{L^6}\\
&~~~ +\|{\cal C}^{n+1}-\Pi_h^n{\cal C}^{n+1}\|_{L^3}
\Big\|\nabla\cdot\Big((D(U^{n+1})-D(U^n))\nabla \psi
\Big)\Big\|_{L^{3/2}}
\\
&\leq C\|D(U^{n+1})-D(U^n)\|_{L^2} h^2\|\psi\|_{W^{2,6}}\\
&~~~+Ch^2\|{\cal C}^{n+1}\|_{W^{2,3}}\Big(\|\nabla
(D(U^{n+1})-D(U^n))\|_{L^{3/2}} \|\nabla
\psi\|_{L^\infty}+\|U^{n+1}-U^n\|_{L^2}\|\psi\|_{W^{2,6}}\Big)\\
&\leq C \left ( \| D_tU^{n+1}\|_{L^2}+C\|\nabla
D_tU^{n+1}\|_{L^{3/2}}+\|\nabla U^{n+1}
D_tU^{n+1}\|_{L^{3/2}} \right ) \tau h^2\|\psi\|_{H^3}\\
&\leq C \left (
\| D_tU^{n+1}\|_{L^2}+C\|\nabla
D_tU^{n+1}\|_{L^{3/2}}+\|\nabla U^{n+1} D_tU^{n+1}\|_{L^{3/2}}
\right ) )
\tau h^2\|\varphi\|_{H^1} .
\end{align*}
Therefore,
\begin{align*}
&\|\Pi_h^n{\cal C}^{n+1}-\Pi_h^{n+1}{\cal C}^{n+1}\|_{H^{-1}}\\
&\leq C \left ( \| D_tU^{n+1}\|_{L^2}+C\|\nabla
D_tU^{n+1}\|_{L^{3/2}}+\|\nabla U^{n+1} D_tU^{n+1}\|_{L^{3/2}}
\right ) ) \tau h^2 .
\end{align*}
By (\ref{errorESTTD}), we have
\begin{align*}
&\| D_tU^{n+1}\|_{L^2}\leq \| D_te_u^{n+1}\|_{L^2}+\|
D_t {\bf u}^{n+1}\|_{L^2}
\leq C  +\| D_t {\bf u}^{n+1}\|_{L^2} ,\\[5pt]
&\|\nabla U^{n+1} D_tU^{n+1}\|_{L^{3/2}}\leq \|\nabla
U^{n+1}\|_{L^6} \|D_tU^{n+1}\|_{L^2}\leq C + C\|
D_t {\bf u}^{n+1}\|_{L^2} .
\end{align*}
With the regularity assumption (\ref{StrongSOlEST}) on ${\bf u}$ and
the estimate (\ref{errorESTTD}), we derive that
\begin{align*}
\biggl(\sum_{n=0}^{N-1}\tau\|\Pi_h^n{\cal C}^{n+1}-\Pi_h^{n+1}{\cal
C}^{n+1}\|_{H^{-1}}^2\biggl)^{\frac{1}{2}}\leq C\tau h^2 .
\end{align*}
Since
\begin{align*}
&\|D_t({\cal C}^{n+1} -\Pi_h^{n+1}{\cal C}^{n+1})\|_{H^{-1}}\\
&\leq \|D_t{\cal C}^{n+1} -\Pi_h^nD_t{\cal C}^{n+1}\|_{H^{-1}}
+ \tau^{-1} \|\Pi_h^{n+1}{\cal C}^{n+1}-\Pi_h^n{\cal C}^{n+1}\|_{H^{-1}}  \\
&\leq C\|D_t{\cal C}^{n+1}\|_{H^2}h^2 + \tau^{-1} \|\Pi_h^{n+1}{\cal
C}^{n+1}-\Pi_h^n{\cal C}^{n+1}\|_{H^{-1}} \\
&\leq C\|D_te_c^{n+1}\|_{H^2}h^2
+C\|D_tc^{n+1}\|_{H^2}h^2 +
\tau^{-1} \|\Pi_h^{n+1}{\cal C}^{n+1}-\Pi_h^n{\cal
C}^{n+1}\|_{H^{-1}},
\end{align*}
with (\ref{StrongSOlESTTD}) and
(\ref{StrongSOlEST}), (\ref{knldfahlae}) follows immediately.

The
proof of Lemma \ref{njekfqhklad} is complete.
\endproof
\bigskip
\begin{theorem}\label{ErrestFEM-TD}
{\it Suppose that the initial-boundary value problem
{\rm (\ref{e-fuel-1})-(\ref{e-fuel-4})} has a unique solution $(p, {\bf
u},c)$ which satisfies {\rm (\ref{StrongSOlEST})}. Then there exist
positive constants $h_0$ and $\tau_0$ such that when $h<h_0$ and
$\tau < \tau_0$, the finite element system
{\rm (\ref{e-FEM-1})-(\ref{e-FEM-2})} admits a unique solution $(P^n_h,
U^n_h , {\cal C}^n_h )$, $n=1,\cdots,N$, which satisfies
\begin{align*}
\max_{1\leq n\leq N}\|P^n_h - P^n\|_{L^2}
+\max_{1\leq n\leq N}\| U^n_h - U^n\|_{L^2}
+\max_{1\leq n\leq N}\|{\cal C}^n_h -{\cal C}^n\|_{L^2} \leq C h^2 .
\end{align*}
where $\{(P^n, U^n , {\cal C}^n )\}_{n=1}^N$ is the solution
of the time-discrete system {\rm (\ref{TDe-fuel-1})-(\ref{TDBC})}
. }
\end{theorem}
\medskip
\noindent{\it Proof}~~~ The solution to the time-discrete system
(\ref{TDe-fuel-1})-(\ref{TDBC}) satisfies
\begin{align}
& \biggl(\frac{\mu({\cal C}^{n})}{k(x)} U^{n+1},\,v_h\biggl)
=-\Big(P^{n+1} ,\, \nabla \cdot v_h \Big),
\label{tde-FEM-1}\\[3pt]
& \Big(\nabla\cdot U^{n+1} ,\, \varphi_h\Big) =\Big(q^I-q^P,\,
\varphi_h\Big),
\label{tde-FEM-2}\\[3pt]
& \Big(\Phi  D_t{\cal C}^{n+1}, \, \phi_h\Big)
+ \Big(D(U^{n+1})\nabla {\cal C}^{n+1}, \, \nabla \phi_h \Big) \nonumber\\
&~~~ + \Big( U^{n+1}\cdot\nabla {\cal C}^{n+1},\, \phi_h\Big)= \Big(\hat
c q^I-{\cal C}^{n+1} q^P, \, \phi_h\Big). \label{tde-FEM-3}
\end{align}
for any $v_h\in H_h$, $\varphi_h\in S_h$ and $\phi_h\in V_h$. The
above equations with the finite element system
(\ref{e-FEM-1})-(\ref{e-FEM-3}) imply that the error functions
$P^{n+1}_h-\Pi_h P^{n+1} $, $U^{n+1}_h- U^{n+1}$, ${\cal
C}^{n+1}_h-{\cal C}^{n+1}$ satisfy
\begin{align}
& \biggl(\frac{\mu({\cal C}^{n}_h)}{k(x)} U^{n+1}_h-\frac{\mu({\cal
C}^{n})}{k(x)} U^{n+1},\,v_h\biggl) =-\Big(P^{n+1}_h-\Pi_hP^{n+1},\,
\nabla \cdot v_h \Big),
\label{erre-FEM-1}\\[3pt]
& \Big(\nabla\cdot( U^{n+1}_h- U^{n+1}) ,\, \varphi_h\Big) =0,
\label{erre-FEM-2}\\[3pt]
& \Big(\Phi D_t({\cal C}^{n+1}_h-{\cal C}^{n+1}), \, \phi_h\Big) +
\Big(D( U^{n+1}_h)\nabla({\cal C}^{n+1}_h-\Pi_h^{n+1}{\cal
C}^{n+1}),
 \, \nabla \phi_h \Big) \nn\\
& = -\Big( U^{n+1}\cdot\nabla ({\cal C}^{n+1}_h-{\cal C}^{n+1}),\,
\phi_h\Big) -\Big(( U^{n+1}_h- U^{n+1})\cdot\nabla {\cal
C}^{n+1}_h,\, \phi_h\Big)
\nn\\
&~~~ -\Big(({\cal C}^{n+1}_h-{\cal C}^{n+1}) q^P, \,
\phi_h\Big)+\Big((D( U^{n+1})-D( U^{n+1}_h))\nabla\Pi_h^{n+1}{\cal
C}^{n+1},
 \, \nabla \phi_h \Big) .
\label{erre-FEM-3}
\end{align}

First, we present an upper bound for the error function $\|
U^{n+1}_h- U^{n+1}\|_{L^2}$ in terms of $\|{\cal C}^{n+1}_h-{\cal
C}^{n+1}\|_{L^2}$. By the definition of the projection operator
$Q_h$ in (\ref{Q}), from (\ref{erre-FEM-2}) we see that
$$
\Big(\nabla\cdot( U^{n+1}_h-Q_h U^{n+1}) ,\, \varphi_h\Big) =0,
\quad\mbox{for all}~~\varphi_h\in S_h,
$$
which implies that
$\nabla\cdot(U^{n+1}_h-Q_h U^{n+1})=0$ in $\Omega$. Taking
$v_h=U^{n+1}_h-Q_h U^{n+1}$ in (\ref{erre-FEM-1}), we get
\begin{align*}
\biggl(\frac{\mu({\cal C}^{n}_h)}{k(x)}\big( U^{n+1}_h-Q_h
U^{n+1}\big) +\frac{\mu({\cal C}^{n}_h)}{k(x)}\big(Q_h U^{n+1}-
U^{n+1}\big)
\nn\\
+\biggl(\frac{\mu({\cal C}^{n}_h)}{k(x)} -\frac{\mu({\cal
C}^{n})}{k(x)}\biggl) U^{n+1},~  U^{n+1}_h -Q_h U^{n+1}\biggl)=0,
\end{align*}
which further implies that
\begin{align*}
\big\| U^{n+1}_h - Q_h U^{n+1}\big\|_{L^2} & \leq C\big\|Q_h
U^{n+1}- U^{n+1}\big\|_{L^2} +C\big\|{\cal C}^{n}_h-{\cal
C}^{n}\big\|_{L^2}.
\end{align*}
With (\ref{p-app}), the above inequality reduces to
\begin{align}\label{errestU0}
\big\| U^{n+1}_h - U^{n+1}\big\|_{L^2} \leq C h^2  + C\big\|{\cal
C}^{n}_h-{\cal C}^{n}\big\|_{L^2} \, .
\end{align}

Secondly, we take $\phi_h= {\cal C}^{n+1}_h-\Pi_h^{n+1}{\cal
C}^{n+1}$ in (\ref{erre-FEM-3}) and
get
\begin{align}
& \frac{1}{2} D_t \biggl \| \sqrt{\Phi}\big({\cal
C}^{n+1}_h-\Pi_h^{n+1} {\cal C }^{n+1}\big) \biggl \|^2_{L^2} +
\biggl \| \sqrt{D(U^{n+1}_h)}\nabla({\cal C}^{n+1}_h-\Pi_h^{n+1}
{\cal C}^{n+1}) \biggl \|^2_{L^2}
\nn\\
& \leq \big(\Phi D_t({\cal C}^{n+1}-\Pi_h^{n+1} {\cal
C}^{n+1}),\phi_h-\mbox{
$\frac{1}{|\Omega|}\int_\Omega\phi_h\d x$}\big)
\nn \\
& ~~~ + C \|q^I - q^P \|_{L^3} \| {\cal C}^{n+1}_h-\Pi_h^{n+1}{\cal
C}^{n+1}\|_{L^6} \| {\cal C}^{n+1}_h-{\cal C}^{n+1} \|_{L^2}
\nn \\
& ~~~ +C \| U^{n+1} \|_{L^\infty} \| \nabla ({\cal
C}^{n+1}_h-\Pi_h^{n+1}{\cal C}^{n+1}) \|_{L^2} \|{\cal
C}^{n+1}_h-{\cal C}^{n+1}\|_{L^2}
\nn\\
&~~~ +C \|{\cal C}^{n+1}_h-\Pi_h^{n+1}{\cal C}^{n+1}\|_{L^6} \big (
\| U^{n+1}_h - U^{n+1} \|_{L^2} \| \nabla (\Pi_h^{n+1}{\cal
C}^{n+1})\|_{L^3}
\nn \\
&~~~+ \| U^{n+1}_h - U^{n+1} \|_{L^2} \|\nabla ({\cal
C}^{n+1}_h-\Pi_h^{n+1}{\cal C}^{n+1} )\|_{L^3} \big)
\nn\\
&~~~+ C\|q^P\|_{L^3}(\|{\cal C}^{n+1}_h-\Pi_h^{n+1}{\cal
C}^{n+1}\|_{L^3}^2+ \|{\cal C}^{n+1}_h-\Pi_h^{n+1}{\cal
C}^{n+1}\|_{L^6}\|{\cal C}^{n+1}-\Pi_h^{n+1}{\cal
C}^{n+1}\|_{L^2}) \nonumber\\
&~~~ + C \| \nabla \Pi_h^{n+1} {\cal C}^{n+1} \|_{L^\infty} \|\nabla
({\cal C}^{n+1}_h-\Pi_h^{n+1}{\cal C}^{n+1})\|_{L^2}
\| U^{n+1}_h - U^{n+1} \|_{L^2}\nonumber\\
& \leq (\epsilon+ C h^{-d/6}\big\| U^{n+1}_h - U^{n+1}\big\|_{L^2})
\|\nabla ({\cal C}^{n+1}_h-\Pi_h^{n+1}{\cal
C}^{n+1})\|_{L^2}^2\nn\\
&~~~+ C\epsilon^{-1}\big \|D_t({\cal C}^{n+1}-\Pi_h^{n+1} {\cal
C}^{n+1}) \big \|_{H^{-1}}^2
+C\epsilon^{-1} \|U^{n+1}_h-U^{n+1}\|_{L^2}^2
\nn \\
& ~~~ + C\epsilon^{-1} \|{\cal C}^{n+1}_h-\Pi_h^{n+1}{\cal
C}^{n+1}\|_{L^2}^2 +  + C \epsilon^{-1}h^4 ,
\label{nigew}
\end{align}
where we have used (\ref{p-app}) and the following
estimate:
\begin{align*}
& \left | \Big( U^{n+1}\cdot\nabla ({\cal C}^{n+1}_h-{\cal
C}^{n+1}),\, {\cal C}^{n+1}_h-\Pi_h^{n+1}{\cal C}^{n+1}\Big) \right
|
\\
&~~~~~ = \left | \Big( (q^I - q^P) ({\cal C}^{n+1}_h-{\cal
C}^{n+1}),\, {\cal C}^{n+1}_h-\Pi_h^{n+1}{\cal C}^{n+1} \Big) \right
|
\nn \\
&~~~~~~~ + \left | \Big( U^{n+1} ({\cal C}^{n+1}_h-{\cal
C}^{n+1}),\, \nabla ({\cal C}^{n+1}_h-\Pi_h^{n+1}{\cal C}^{n+1})
\Big) \right |
\\
&~~~~~ \le C \|q^I - q^P \|_{L^3} \| {\cal
C}^{n+1}_h-\Pi_h^{n+1}{\cal C}^{n+1} \|_{L^6} \| {\cal
C}^{n+1}_h-{\cal C}^{n+1} \|_{L^2}
\\
&~~~~~~~ + C \| U^{n+1} \|_{L^\infty} \| \nabla ({\cal
C}^{n+1}_h-\Pi_h^{n+1}{\cal C}^{n+1}) \|_{L^2} \| {\cal
C}^{n+1}_h-{\cal C}^{n+1} \|_{L^2} \, .
\end{align*}
From (\ref{errestU0}) we observe that (\ref{nigew}) reduces to
\begin{align*}
& D_t  \biggl \| \sqrt{\Phi}( {\cal C}^{n+1}_h-\Pi_h^{n+1}{\cal
C}^{n+1}) \biggl \|_{L^2}^2  +  \frac{1}{2}\biggl \|\sqrt{D(U^{n+1}_h)}
\nabla({\cal C}^{n+1}_h-\Pi_h^{n+1}{\cal C}^{n+1}) \biggl \|^2_{L^2}
\\
& \leq Ch^4  +C h^{-1/2} \|{\cal C}^{n}_h-\Pi_h^{n}{\cal C}^{n}\|_{L^2}
\|\nabla ({\cal C}^{n+1}-\Pi_h^{n+1}{\cal C}^{n+1})\|_{L^2}^2
\\
&+C(\|{\cal C}^{n+1}_h-\Pi_h^{n+1}{\cal C}^{n+1}\|_{L^2}^2+\|{\cal
C}^{n}_h-\Pi_h^{n}{\cal C}^{n}\|_{L^2}^2) +C\big\|D_t({\cal
C}^{n+1}-\Pi_h^{n+1} {\cal C}^{n+1}) \big \|_{H^{-1}}^2 .
\end{align*}
By applying Gronwall's inequality (with mathematical induction on $\|{\cal C}^{n}_h-\Pi_h^{n}{\cal C}^{n}\|_{L^2}\leq h$), we deduce that
there exists a positive constant $h_0$ such that when $h<h_0$ we have
\begin{align}
\|{\cal C}^{n+1}_h-\Pi_h^{n+1}{\cal C}^{n+1}\|_{L^2} \leq Ch^2 .
\label{nigew3}
\end{align}
From (\ref{p-app}) and (\ref{errestU0}), we further get
\begin{align}
&\big\| U^{n+1}_h- U^{n+1}\big\|_{L^2} \leq Ch^2,\\
&\big\|{\cal C}^{n+1}_h-{\cal C}^{n+1}\big\|_{L^2}\leq Ch^2 .
\label{p-c}
\end{align}

Finally, we estimate the error $\| P_h-P^{n+1}\|_{L^2}$. We redefine
$g$ to be  the solution to the equation
\begin{align*}
&-\nabla\cdot\biggl(\frac{k(x)}{\mu({\cal C}^{n})}\nabla g\biggl)
=P^{n+1}_h-\Pi_hP^{n+1}
\end{align*}
with the boundary condition $\frac{k(x)}{\mu({\cal C}^{n})}\nabla
g\cdot {\bf n}=0$ on $\partial\Omega$. Easy to check that
$$\|g\|_{H^2}\leq C\|P^{n+1}_h-\Pi_hP^{n+1}\|_{L^2}.$$ Let
$$
v_h=Q_h\biggl(\frac{k(x)}{\mu({\cal C}^{n})}\nabla g\biggl)
$$
Then
$$
(\varphi,\nabla \cdot v_h)  =-
(\varphi,P^{n+1}_h-\Pi_hP^{n+1}),~~\varphi\in S_h
$$
and from (\ref{erre-FEM-1}) we obtain
\begin{align*}
\|P^{n+1}_h-\Pi_hP^{n+1}\|_{L^2}^2 & = \biggl(\frac{\mu({\cal
C}^{n}_h)}{k(x)} U^{n+1}_h -\frac{\mu({\cal C}^{n})}{k(x)} U^{n+1},
Q_h\biggl(\frac{k(x)}{\mu({\cal C}^{n})}\nabla g\biggl) \biggl)
\nn \\
& \leq C(\|{\cal C}^n_h-{\cal C}^n\|_{L^2}
+\|U^{n+1}_h-U^{n+1}\|_{L^2})
\biggl\|\frac{k(x)}{
\mu({\cal C}^{n})}\nabla g\biggl\|_{H^1}\\
&\leq Ch^2\|P^{n+1}_h-\Pi_hP^{n+1}\|_{L^2} ,
\end{align*}
which implies that
\begin{align*}
\|P^{n+1}_h- \Pi_h P^{n+1}\|_{L^2} \leq Ch^2 .
\end{align*}
The proof of Theorem \ref{ErrestFEM-TD} is complete.
\endproof

\medskip
Combining Theorems \ref{ErrestTDSol} and Theorem \ref{ErrestFEM-TD},
we complete the proof of Theorem \ref{ErrestFEMSol}. \endproof
\section{Numerical examples}
\setcounter{equation}{0}
In this section, we present two numerical examples to confirm our
theoretical analysis.
All computations are performed by using the software FreeFEM++.\medskip

\noindent{\it Example 5.1}~~
We rewrite the system (\ref{e-fuel-1})-(\ref{e-fuel-4})
by
\begin{align}
&\frac{\partial c}{\partial t}-\nabla\cdot(D({\bf u})\nabla
c)+{\bf u}\cdot\nabla c= g,
\label{dfnhsdfui001}\\[3pt]
&\nabla\cdot\u=f,
\label{dfnhsdfui002}\\[3pt]
&\u=-\frac{1}{\mu(c)}\nabla p,
\label{dfnhsdfui003}
\end{align}
where $\Omega = (0,1)\times(0,1)$ and $D({\bf u})=1+|{\bf u}|^2/(1+|{\bf
u}|)$
and $\mu(c)=1+c^2$.
The functions $f$ and $g$
are chosen corresponding to the exact solution
\begin{align}
&p=1+1000x^2(1-x)^3y^2(1-y)^3t^2e^{t},
\label{dfnhs001}\\
&{\bf u}=-\frac{1}{\mu(c)}\nabla p,
\label{dfnhs002}\\
&c=0.1+50x^2(1-x)^2y^2(1-y)^2te^{t} .
\label{dfnhs003}
\end{align}
Clearly, the boundary condition (\ref{e-fuel-4}) is satisfied.

A uniform triangular partition with $M+1$ nodes in each direction
is used to generate the FEM mesh (with $h=1/M$). We solve the system by the
proposed method up to the time $t=1$.
To illustrate our error estimates, numerical errors with $\tau=8h^2$
are presented in Table \ref{linear-L2-1}, from which we can see that the $L^2$ errors
are proportional to $O(h^2)$.
To demonstrate the unconditionaly convergence (stability) of the numerical method,
we solve the system with a fixed $\tau$ and several different spatial mesh size $h$.
We present numerical errors in Table \ref{linear-L2-2}. We can observe from
Table \ref{linear-L2-2} that numerical errors behave like $O(\tau)$
as $h/\tau \rightarrow 0$.
This implies that the time-step conditions are not necessary.
\vskip0.1in

\begin{table}[p]
\vskip-0.2in
\begin{center}
\caption{Errors of the Galerkin-mixed FEM in $L^2$ norm.}\vskip 0.1in
\label{linear-L2-1}
\begin{tabular}{l|l|c|ccc}
\hline
$\tau$  &  $h$
& $\| U_h^N - {\bf u}(\cdot ,t_N) \|_{L^2}$
& $\| {\cal C}_h^N - c(\cdot ,t_N) \|_{L^2}$        \\
\hline
1/8   & 1/8
& 2.024E-01  & 7.114E-02 \\
1/32  & 1/16
& 5.264E-02  & 1.713E-02 \\
1/128 & 1/32
& 1.333E-02  & 4.070E-03 \\
\hline
\multicolumn{2}{c|}{convergence rate}
& 1.98       & 2.07 \\
\hline
\end{tabular}
\end{center}

\begin{center}
\vskip-0.1in
\caption{Errors of the Galerkin-mixed FEM with fixed $\tau$ and refined $h$.
}\vskip 0.1in
\label{linear-L2-2}
\begin{tabular}{l|c|c|ccc}
\hline
$\tau=0.05$ & $h$
& $\| U_h^N - {\bf u}(\cdot ,t_N) \|_{L^2}$ & $\| {\cal C}_h^N - c(\cdot
,t_N) \|_{L^2}$         \\
\cline{2-4}
 & 1/8  
 & 1.955E-01  & 4.748E-02 \\
 & 1/16 
 & 5.531E-02  & 2.081E-02 \\
 & 1/32 
 & 2.409E-02  & 1.077E-02 \\
 & 1/64 
 & 1.998E-02  & 8.243E-03 \\
\hline
\hline
$\tau=0.1$  &  $h$   
 & $\| U_h^N - {\bf u}(\cdot ,t_N) \|_{L^2}$ & $\| {\cal C}_h^N - c(\cdot
,t_N) \|_{L^2}$         \\
\cline{2-4}
 & 1/8   
 & 1.998E-01  & 6.216E-02 \\
 & 1/16  
 & 6.577E-02  & 3.348E-02 \\
 & 1/32  
 & 4.168E-02  & 2.240E-02 \\
 & 1/64  
 & 3.910E-02  & 1.961E-02 \\
\hline
\hline
$\tau=0.25$ &  $h$ 
  & $\| U_h^N - {\bf u}(\cdot ,t_N) \|_{L^2}$ & $\| {\cal C}_h^N - c(\cdot
,t_N) \|_{L^2}$         \\
\cline{2-4}
 & 1/8   
 & 2.195E-01  & 1.336E-01 \\
 & 1/16  
 & 1.088E-01  & 9.885E-02 \\
 & 1/32  
 & 9.491E-02  & 8.426E-02 \\
 & 1/64  
 & 9.349E-02  & 8.062E-02 \\
\hline
\end{tabular}
\end{center}

\begin{center}
\caption{Errors of the Galerkin-mixed FEM with fixed $\tau$ and refined
$h$.}\vskip 0.1in
\label{linear-L2-5}
\begin{tabular}{l|c|c|ccc}
\hline
$\tau=0.05$ & $M$ 
  & $\| U_h^N - {\bf u}(\cdot ,t_N) \|_{L^2}$ & $\| {\cal C}_h^N - c(\cdot
,t_N) \|_{L^2}$         \\
\cline{2-4}
 & $32$  
 & 7.105E-02  & 1.445E-02 \\
 & $64$  
 & 2.526E-02  & 4.022E-03 \\
 & $128$  
 & 1.523E-02  & 7.754E-04 \\
\hline
\hline
$\tau=0.1$ & $M$ 
  & $\| U_h^N - {\bf u}(\cdot ,t_N) \|_{L^2}$ & $\| {\cal C}_h^N - c(\cdot
,t_N) \|_{L^2}$         \\
\cline{2-4}
 & $32$  
 & 7.560E-02  & 1.569E-02 \\
 & $64$  
 & 3.523E-02  & 4.340E-03 \\
 & $128$  
 & 2.869E-02  & 1.248E-03 \\
\hline
\hline
$\tau=0.25$ & $M$ 
  & $\| U_h^N - {\bf u}(\cdot ,t_N) \|_{L^2}$ & $\| {\cal C}_h^N - c(\cdot
,t_N) \|_{L^2}$         \\
\cline{2-4}
 & $32$  
 & 9.719E-02  & 2.900E-02 \\
 & $64$ 
 & 6.960E-02  & 1.429E-02 \\
 & $128$ 
 & 6.632E-02  & 7.940E-03 \\
\hline
\end{tabular}
\end{center}
\end{table}

\begin{figure}[h]
\begin{minipage}[b]{0.38\linewidth}
\centering
\includegraphics[width=\textwidth]{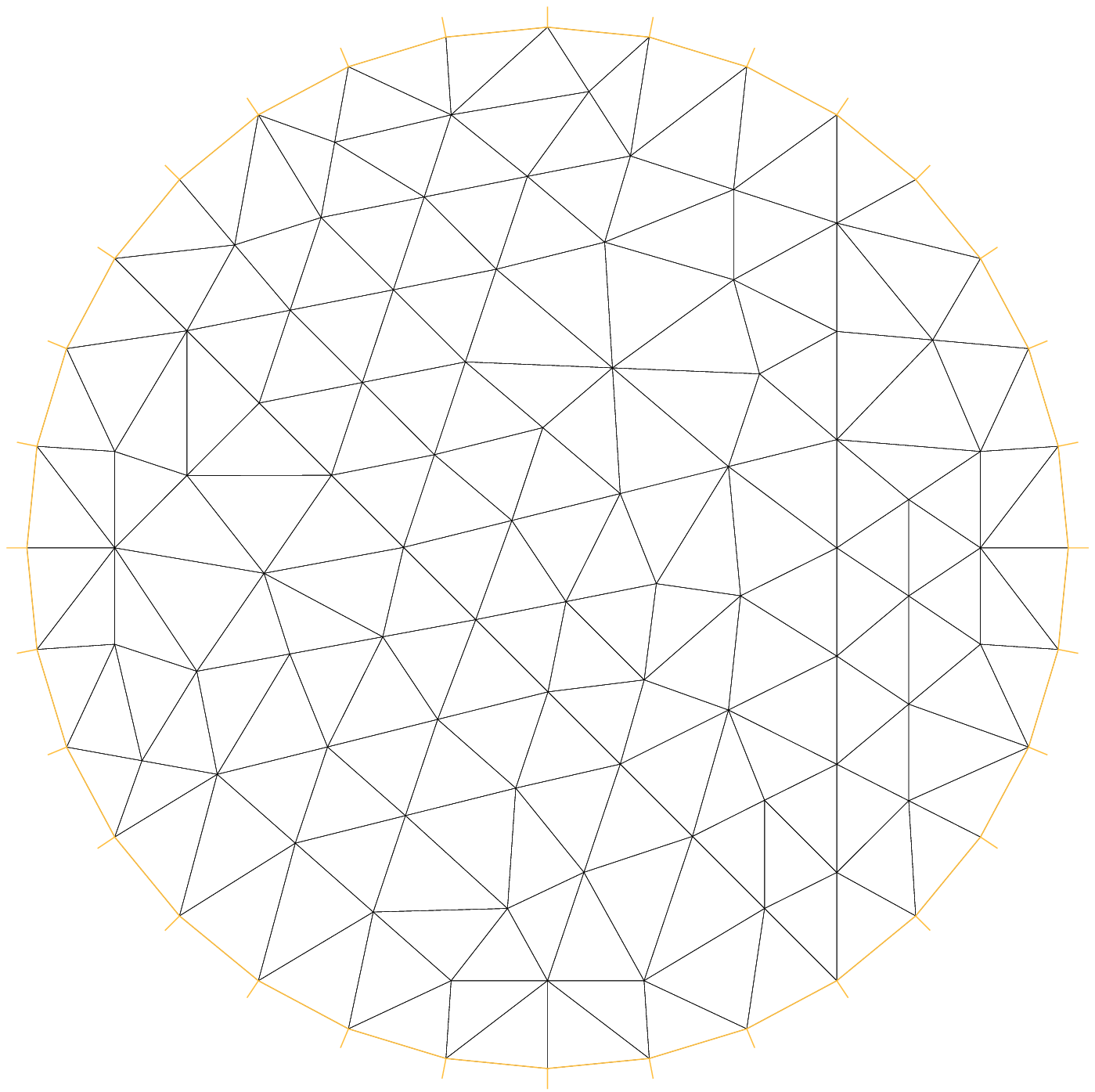}
\end{minipage}
\hspace{-1.8cm}
\begin{minipage}[b]{0.38\linewidth}
\centering
\includegraphics[width=\textwidth]{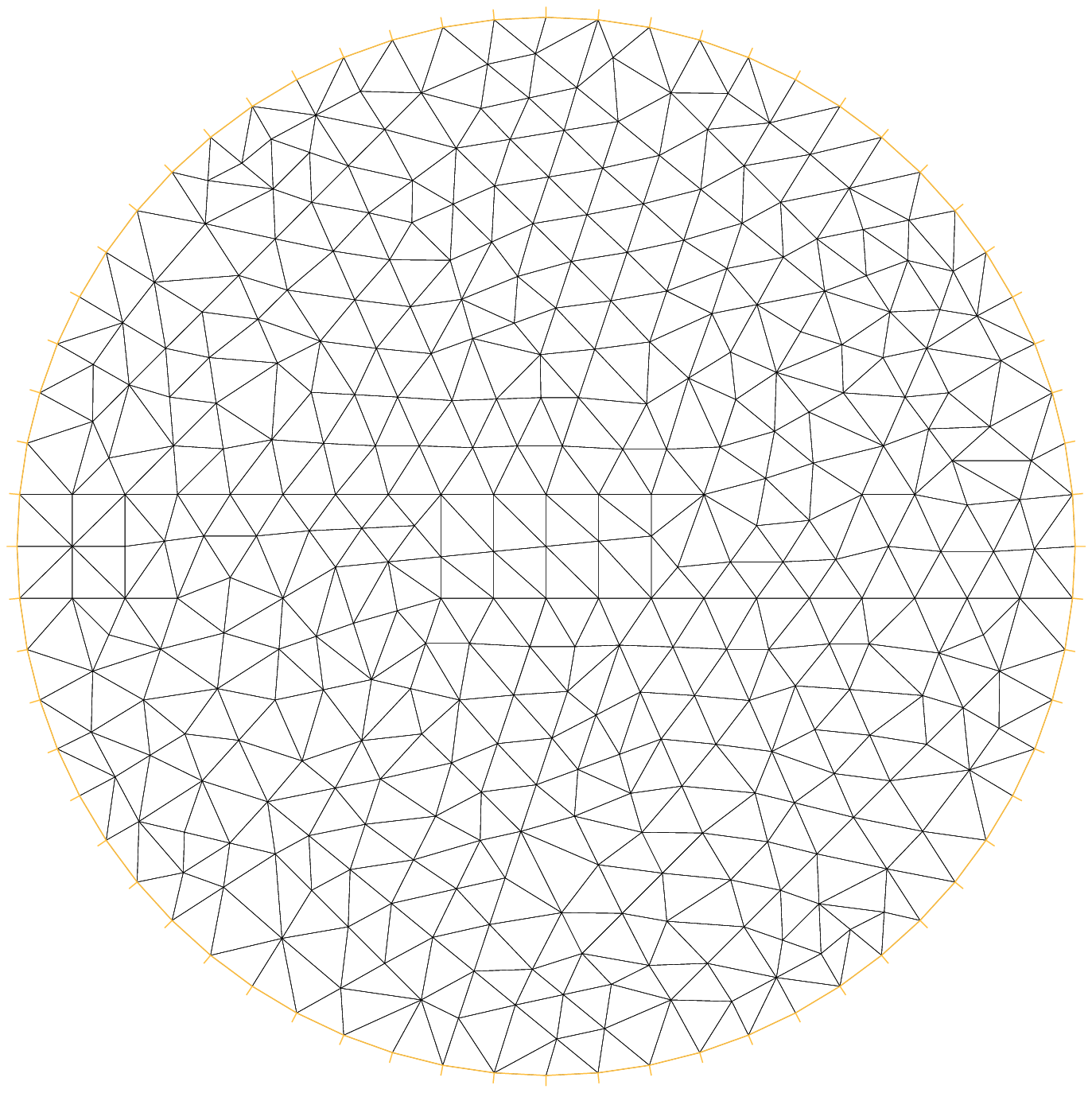}
\end{minipage}
\hspace{-1.8cm}
\begin{minipage}[b]{0.38\linewidth}
\centering
\includegraphics[width=\textwidth]{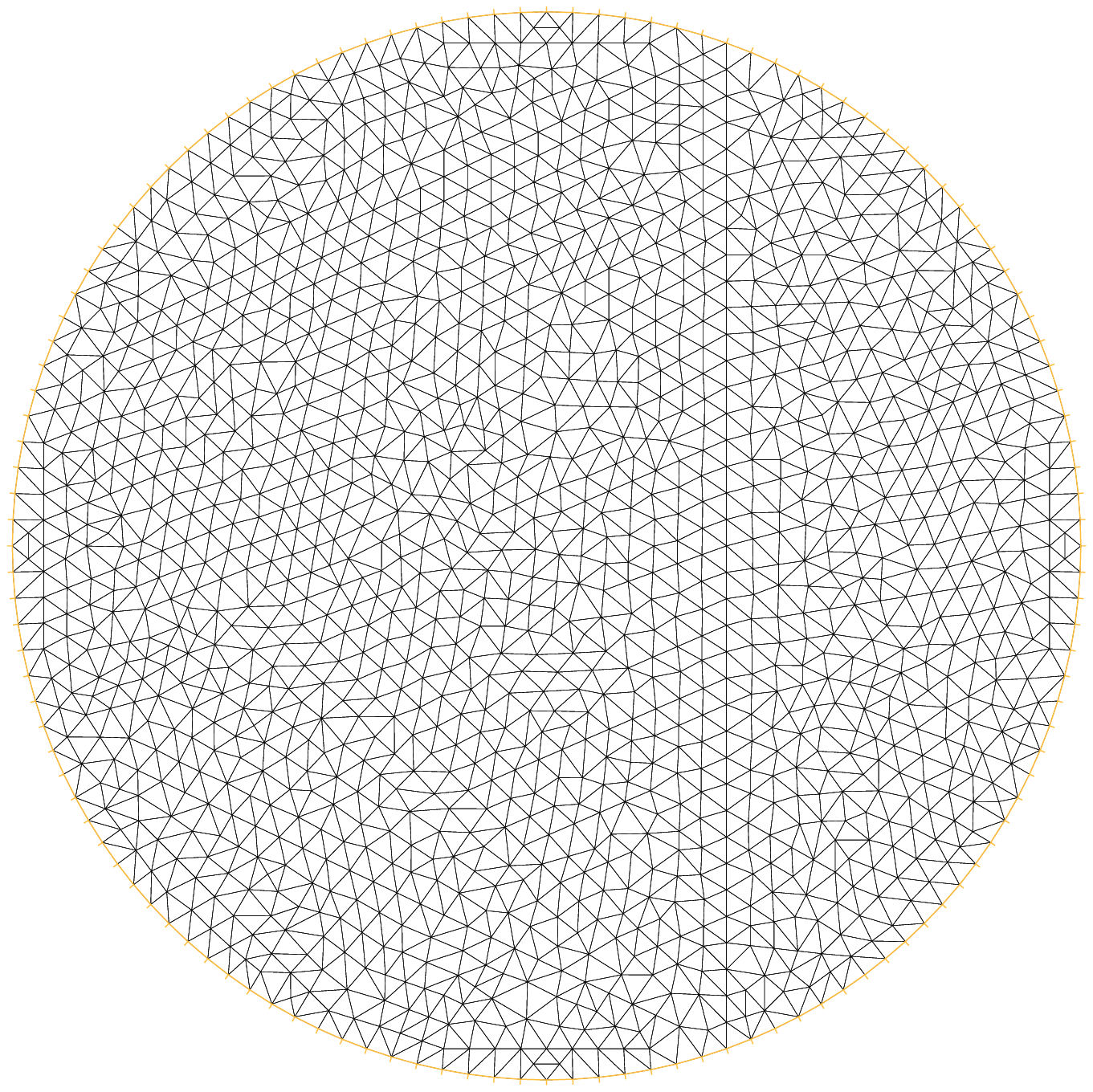}
\end{minipage}
\caption{The FEM meshes with $M=32$, $M=64$ and $M=128$.}\label{figure1}
\end{figure}

\noindent{\it Example 5.2}~~ We consider the
equations (\ref{dfnhsdfui001})-(\ref{dfnhsdfui003})
in a circle centered at $(0.5,0.5)$ with the radius
$0.5$ and with inhomogeneous Neumann boundary
conditions correspondingly to the exact solution
given in (\ref{dfnhs001})-(\ref{dfnhs003}).
The mesh generated here consists of $M$
boundary points with $M=32,64,128$, respectively, as shown in Figure
\ref{figure1}.
Numerical errors with fixed $\tau$ and several different $h$
are presented in Table \ref{linear-L2-5}, which
also show clearly that no time-step condition is needed.
\vskip0.1in

\section{Conclusions}
We have studied error analysis for a nonlinear and strongly coupled
parabolic system from incompressible miscible flow in porous media
with a commonly-used Galerkin-mixed FEM and linearized semi-implicit
Euler scheme. Optimal $L^2$ error estimates were obtained almost
without any time-step condition, while all previous
works imposed certain restriction for the time-step size.
The unconditional error analysis presented in this paper
can be extended to models with other boundary conditions and numerical
methods
with high-order approximations, while here we only focus our analysis on
the problem with a homogeneous boundary condition and the lowest-order
Galerkin-mixed FEM.
In fact, we have proved the error estimates:
$$
\|U^n_h-U^n\|_{L^2}\leq Ch^2,\qquad \|U^n-u^n\|_{L^2}\leq C\tau \, .
$$
We can see from our proof that the two inequalities also hold
for higher-order finite element methods.
The inequalities imply the boundedness of numerical solution.
With some more precise analysis for the time-discrete system,
optimal error estimates of high-order Galerkin type methods
can be obtained in the traditional way.
Also we believe
that the idea of the error splitting coupled with the regularity
analysis of the time-discrete PDEs can be applied to many other
nonlinear parabolic PDEs and time discretizations to obtain optimal
error estimates unconditionally.

\bigskip

\noindent{\bf Acknowledgements}~ The authors would like to
thank the anonymous referees for their valuable suggestions and comments.

\end{document}